\documentstyle[12pt,leqno,amscd,amssymb,pstricks,latexsym,amsbsy,xypic,mathrsfs,verbatim]{amsart}
\numberwithin{equation}{section}

\theoremstyle{plain}
\newtheorem{prop}{Proposition}[section]
\newtheorem{thm}[prop]{Theorem}
\newtheorem{cor}[prop]{Corollary}
\newtheorem{lem}[prop]{Lemma}
\newtheorem{defn}[prop]{Definition}

\theoremstyle{definition}
\newtheorem{example}[prop]{Example}
\newtheorem{rem}[prop]{Remark}

\newcommand{\ff}{\footnote}

\def\rad{{\rm rad}}

\def\t1{{S_\vartriangle(2,r)}}

\def\sn{{S_\vartriangle(n,r)}}
\def\bbz{{\mathbb{Z}}}\def\bbq{{\mathbb{Q}}}
\def\bbn{{\mathbb{N}}}\def\bbc{{\mathbb{C}}}
\def\lb{{\lambda}}
\def\c2{{\Theta_\vartriangle(2,2)}}
\def\fkS{{\frak S}}
\def\fsr{{\frak S_{\vartriangle}}}

\def\c{{e_C}}

\def\lr{{\Lambda_\vartriangle(n,r)}}

\def\diag{{\rm diag}}

\def\h{H_0(r)}

\def\bbf{{\mathbb F}}
\def\wfkF{\fkF_\vartriangle}

\def\fkF{{\frak F}}
\def\udim{{\underline {\rm dim}}\,}

\def\wfkFn{\fkF_{\vartriangle, n}}

\def\row{{\rm row}} \def\col{{\rm col}}

\def\cnr{{\Theta_\vartriangle(n,r)}}

\def\ggp#1#2{\left[\kern-3.2pt\left[{#1\atop #2}\right]\kern-3.2pt\right]}

\def\uj{\underline{j}}

\def\ui{\underline{i}}\def\ul{\underline{l}}

\def\c{\mathcal{C}}

\def\h2{H_\vartriangle(2)}

\def\c{\mathbb{C}}

\def\flb{{\frak l}_{\lb}}

\def\fl{{\frak l}_{\lb^{(1)}}}

\def\bfli{\bar{\frak l}_{\lb^{(i)}}}

\def\bfl{\bar{\frak l}_{\lb}}

\def\bs{\bar{S}}

\def\End{\mbox{\rm End}}

\def\GL{\mbox{\rm GL}}

\def\add{\mbox{\rm add}}

\def\fii{{{\frak S}_{\vartriangle,\underline{i}}}}

\def\fiii{{{\frak S}_{\underline{i}}}}
\def\fjf{{{\frak S}_{\underline{j}}}}

\def\fj{{{\frak S}_{\vartriangle,\underline{j}}}}
\def\fjl{{{\frak S}_{\vartriangle,\underline{j},\underline{l}}}}

\def\fjlf{{{\frak S}_{\underline{j},\underline{l},\varepsilon'}}}

\def\fijf{{{\frak S}_{\underline{i},\underline{j},\varepsilon}}}
\def\fijff{{{\frak S}_{\underline{i},\underline{j}}}}

\def\fij{{{\frak S}_{\vartriangle,\underline{i},\underline{j}}}}

\def\fil{{{\frak S}_{\vartriangle, \underline{i}, \underline{l}\delta}}}
\def\fijl{{{\frak S}_{\vartriangle, \underline{i},\underline{j}, \underline{l}\delta}}}
\def\ff{{{\frak S}_{\vartriangle, \underline{i}\delta, \underline{l}}}}
\def\fg{{{\frak S}_{\vartriangle, \underline{i}\delta,\underline{j}, \underline{l}}}}

\def\ffj{{{\frak S}_{\underline{j}}}}
\def\ffjl{{{\frak S}_{\underline{j},\underline{l},\varepsilon'}}}
\def\fflj{{{\frak S}_{\underline{i},\underline{j},\varepsilon}}}

\def\ww{{{\frak S}_{\underline{i},\underline{l}\delta,\varepsilon'\delta+\varepsilon}}}
\def\vv{{{\frak S}_{\underline{i},\underline{j},\underline{l}\delta,\varepsilon'\delta,\varepsilon}}}

\def\wv{{{\frak S}_{\underline{i}\delta,\underline{l},\varepsilon'+\varepsilon\delta}}}
\def\vw{{{\frak S}_{\underline{i}\delta,\underline{j},\underline{l},\varepsilon',\varepsilon\delta}}}

\def\xij{{ \xi_{\underline{i}, \underline{j}+n\varepsilon}}}
\def\xjl{{ \xi_{\underline{j}, \underline{l}+n\varepsilon'}}}

\def\xil{{ \xi_{\underline{i}, \underline{l}\delta+n(\varepsilon'\delta+\varepsilon)}}}

\def\yil{{ \xi_{\underline{i}\delta, \underline{l}+n(\varepsilon'+\varepsilon\delta)}}}

\def\zv{\bbz[v,v^{-1}]}

\def\whS{{\widehat S}}

\begin{document}

\title{Affine quasi-heredity of affine Schur algebras}

\author{ Bangming Deng  and  Guiyu Yang}


\address{ Mathematical Sciences Center, Tsinghua University,
Beijing 100084,  China.} \email{ bmdeng@@math.tsinghua.edu.cn }
\address{School of Science, Shandong University of Technology, Zibo
255049, China} \email{yanggy@@mail.bnu.edu.cn}

\thanks{Supported partially by the Natural Science Foundation of China.}

\date{\today}

\keywords{Affine Schur algebras, Affine quasi-hereditary algebras}

\subjclass[2000]{20G43, 17B37, 20G05}

\begin{abstract}

In this paper we prove that the affine Schur algebra $\whS(n,r)$ is
affine quasi-hereditary. This result is then applied to show that $\whS(n,r)$
has finite global dimension and its centralizer subquotient algebras
are Laurent polynomial algebras. We also use the result to give a parameter set of simple
$\whS(n,r)$-modules and identify this parameter set with that given in
\cite{DDF}.

\end{abstract}

\maketitle

\section{Introduction}

Affine (quantum) Schur algebras play a central role in linking the
representations of affine quantum groups and affine Hecke
algebras. These algebras can be defined in several equivalent ways and has been
widely studied; see, for example, \cite{GV, rmg,L,vv,DDF}.

The notion of affine quasi-hereditary algebras, which are an affine analogue
of quasi-hereditary algebras, was introduced in \cite{KX}
as a description of a kind of algebras which have affine Hecke
algebras of type $A$ as a perfect model. In \cite{K}, affine quasi-hereditary (graded)
algebras was systematically studied and an affine analogue of the
Cline--Parshall--Scott Theorem was given. In particular, quiver Hecke algebras (or KLR algebras) of finite
type and Kato's geometric extension algebras are shown to be affine quasi-hereditary algebras in \cite{K}.

It was shown in \cite{C} that affine quantum Schur algebra $\whS_v(n,r)$
are affine cellular in the sense of K$\ddot{\rm o}$nig--Xi \cite{KX} by proving
that the affine cell ideals of $\whS_v(n,r)$ are generated by canonical bases introduced by Lusztig in
\cite{L} and can be viewed as generalized matrix algebras over the
representation ring of certain direct product of general linear
groups. Recently, Nakajima \cite{N} proved that the cell ideals of
modified quantum affine algebras are idempotent, and, in particular, the
BLN algebras are affine cellular with idempotent cell ideals. As a
kind of BLN algebra, $\whS_v(n,r)$ has idempotent cell ideals
in case $n>r$.

The aim of this paper is to prove that the affine Schur algebra
$\whS(n,r)$ (by evaluating the parameter $v=1$) is affine
quasi-hereditary in the sense of \cite{KX}. In case
$n>r$, we identify the cell ideals given in \cite{C} with the
ideals generated by some idempotents in $\whS(n,r)$. This identification
gives the affine quasi-heredity of $\whS(n,r)$. In case $n\leq r$, we use Schur functors to
derive an affine-heredity chain of $\whS(n,r)$, which affords an affine
quasi-hereditary structure for $\whS(n,r)$ in this case. As applications, we prove
that $\whS(n,r)$ has finite global dimension and its centralizer subquotient algebras of $\whS(n,r)$ are Laurant
polynomial algebras. Furthermore, we use the affine-heredity chain to give a parameter set of simple
$\whS(n,r)$-modules and identify this parameter set with that given in
\cite{DDF}.

In a forthcoming paper, we will consider the affine quasi-heredity of affine quantum Schur algebras.

The paper is organized as follows. In Sections 2--4 we recall the definitions of affine cellular algebras and
affine quasi-hereditary algebras and review two definitions of the affine Schur algebra $\whS(n,r)$ and
a construction of its affine cell chain.
The affine quasi-heredity of $\whS(n,r)$ in case $n>r$ is proved in Section 5. We then use the Schur functor to prove the
affine quasi-heredity of $\whS(n,r)$ in case $n\leq r$ in Section 6. As an application, we prove in Section 7 that the
centralizer subquotient algebras of $\whS(n,r)$ are Laurent polynomial algebras. In the final section, we give
a parameter set of simple $\whS(n,r)$-modules and identify it with the parameter set given in \cite{DDF}.

\section{ Affine cellular algebras and affine quasi-hereditary algebras }

In this section we recall the definitions of affine cellular
algebras and affine quasi-hereditary algebras introduced in
\cite{KX} and \cite{K}, respectively. For a module over an
algebra, we always mean a left module unless otherwise
stated.

Suppose that $k$ is a Noetherian domain. A commutative $k$-algebra
$B$ is called {\it affine} if it is a quotient of a polynomial ring
$k[x_1,\ldots, x_s]$ in finitely many variables $x_1, \ldots, x_s$
by an ideal $I$, i.e. $B=k[x_1,\ldots, x_s]/I$. A $k$-linear
anti-automorphism $\tau$ of a $k$-algebra $A$ with $\tau^2={\rm id}_A$
will be called a $k$-involution on $A$.

\begin{defn} [{\cite{KX}}] \label{def-aff}
Let $A$ be a unitary $k$-algebra with a $k$-linear involution
$\tau$ on $A$. A two-sided ideal $J$ in $A$ is called an affine
cell ideal if the following conditions are satisfied:

 $(1)$ $\tau(J)=J$.

$(2)$ There is a free $k$-module $V$ of finite rank and an affine
algebra $B$ with a $k$-involution $\sigma$ such that
$\Delta=V\otimes_kB$ is an $A$-$B$-bimodule, where the right
$B$-module structure is induced by that of the regular right
$B$-module $B_B$ .

$(3)$ There is an $A$-$A$-bimodule isomorphism
$\alpha:J\rightarrow \Delta\otimes_B\Delta'$, where
$\Delta'=B\otimes_kV$ is a
 $B$-$A$-bimodule with the left $B$-structure induced by $_BB$ and the right $A$-module structure is
given as $(b\otimes v)a=p(\tau(a)(v\otimes b)), $ where $p$ is the
switch map:
$$p:~~\Delta\otimes_B\Delta'\longrightarrow
\Delta'\otimes_B\Delta,$$$$ x\otimes y\longrightarrow y\otimes x,
~~\text{for}~~ x\in\Delta~~\text{and}~~ y\in \Delta'.$$

$(4)$ There is the following commutative diagram:

$$\xymatrix{ J \ar[d]^-{\tau}\ar[r]^-{\alpha}& \Delta \otimes_B \Delta'\ar[d]^-{
v_1\otimes b_1\otimes_B b_2\otimes v_2\rightarrow v_2\otimes \sigma(b_2)\otimes_B \sigma(b_1)\otimes v_1}\\
J\ar[r]^-{\alpha} &\Delta \otimes_B \Delta' }$$

The algebra $A$ is called affine cellular if and only if there is
a $k$-module decomposition $A=J_1'\oplus J_2'\oplus\cdots J_n'$
with $\tau(J_j')=J_j'$ for each $j$, and $J_i=\oplus_{1\leq l\leq
i}J_l'$ gives a chain of two-sided ideals of $A$: $0=J_0\subset
J_1\subset J_2\subset\cdots\subset J_n=A$, and for each $1\leq
i\leq n$, $J_i'=J_i/J_{i-1}$ is an affine cell ideal of
$A/J_{i-1}$.

We call this chain an affine cell chain for the affine cellular
algebra $A$. The bimodule $\Delta$ will be called a cell lattice for the affine
cell ideal $J$.

\end{defn}

The following two lemmas are taken from \cite[Lem. 2.4] {KX} and \cite[Th.~4.3]{KX}.

\begin{lem}\label{ring}

Suppose that $K$ is another Noetherian doamin and $\psi: k\rightarrow K$ is a homomorphism of
rings with identity. If $A$ is an affine cellular $k$-algebra, then $K\otimes_k A$ is an affine cellular
$K$-algebra with the canonical affine cell chain induced from that of $k$-algebra $A$.

\end{lem}

\begin{lem} \label{usu}
Let $J=V\otimes_k B\otimes_kV$ be an idempotent affine cell ideal
in a $k$-algebra $A$ with the cell lattice $\Delta=V\otimes_k B$
as in Definition \ref{def-aff}.

$(1)$ If there is a nonzero idempotent $e$ in $J$, then $_AJ$ is a
projective $A$-module, $J=AeA$, and there is an equality
$\add(_A Ae)=\add(_A \Delta)$.

$(2)$ If $\rad(B)=0$, then $\End(_A \Delta)\cong B$.

\end{lem}

Now we give the definition of an affine quasi-hereditary algebra following \cite{KX,K}.

\begin{defn} \label{quasi}
Let $A$ be a left Noetherian $k$-algebra. An ideal $J$ in $A$
is called an affine-heredity ideal if the following
conditions are satisfied:

\begin{itemize}
\item[(1)] There is an idempotent $e\in J$ such that $J=AeA$;

\item[(2)] As a left $A$-module, $J\cong\underbrace{P\oplus\cdots\oplus
P}_m$, for some projective $A$-module $P$ and some $m\in\bbn$, and
$B :=\End_A(P)$ is an affine algebra, which means that $B$ is a
finitely generated commutative algebra.

\item[(3)] As a right $B$-module, $P$ is finitely generated and flat.
\end{itemize}

The algebra $A$ is called affine quasi-hereditary if there is a
finite chain of ideals
$$0=J_0\subseteq J_1\subseteq\ldots \subseteq J_t=A$$
such that $J_i/J_{i-1}$ is an affine-heredity ideal of $A/J_{i-1}$ for
each $1\leq i\leq t$. Such a chain of ideals is called an
affine-heredity chain.
\end{defn}

\begin{rem}\label{equidef} The definition for an affine quasi-hereditary given here follows from
the one in \cite{KX}. A graded version is given in \cite{K}.
\end{rem}

\section{Affine Schur algebras}

This section is devoted to introducing the geometric definition of affine (quantum)
Schur algebras given by \cite{GV,L}. An algebraic definition of affine Schur algebras
given in \cite{Ya1} will also be reviewed. Finally, we give a
correspondence of the basis elements between these two definitions
and recall some multiplication formulas in \cite{DDF,L,Ya1}.

Let $\bbf$ be a field and let $\bbf[x,x^{-1}]$ be the Laurent polynomial ring in indeterminate
$x$. Fix an $\bbf[x,x^{-1}]$-free module $V$ of rank $r\geq 1$. A lattice in $V$ is a free
$\bbf[x]$-submodule $L$ of $V$ such that $V=L\otimes_{\bbf[x]}\bbf[x,x^{-1}]$.

Let $\wfkF=\wfkFn$ be the set of all cyclic flags
$L=(L_i)_{i\in\bbz}$ of lattices, where each $L_i$ is a lattice in
$V$ such that $L_{i-1}\subseteq L_i$ and $L_{i-n}=x L_i$ for all
$i\in \bbz.$ The group $G$ of automorphisms of the the
$\bbf[x,x^{-1}]$-module $V$ acts on $\wfkF$ by $g\cdot
L=(g(L_i))_{i\in\bbz}$ for $g\in G$ and $L\in \wfkF$. Thus, the
map
\begin{equation}\label{flag} \phi:\wfkF\longrightarrow \lr,\;\; L\longmapsto
\udim L=(\dim_\bbf L_i/L_{i-1})_{i\in \bbz}\end{equation} induces
a bijection between the set of $G$-orbits in $\wfkF$ and $\lr$,
where
$$\lr:=\{(\lb_i)_{i\in\bbz}\mid \lb_i\in\bbn,
\sum_{i=1}^n\lb_i=r \;\mbox{and}\;\lb_i=\lb_{i-n}\; \mbox{for}\;
i\in \bbz\}.$$ We denote by $\frak F_\lb$ the fiber of $\lb$ under
this map $\phi$, i.e. $\frak F_\lb=\phi^{-1}(\lb)$.

Let $$\Lambda(n,r):=\{(\lb_1,\ldots,\lb_n)\mid \lb_i\in\bbn,~
\sum_{1\leq i\leq n}\lb_i=r\}.$$  We usually identify
$\Lambda(n,r)$ with $\lr$ via the following bijection:
$$b: \lr\longrightarrow \Lambda(n,r),~~\lb\longmapsto (\lb_1,\ldots,\lb_n).$$

 The group $G$ also acts diagonally on $\wfkF\times \wfkF$ by $g(L,
L')=(gL, gL')$, where $g\in G$ and $L, L'\in \wfkF$. By
\cite[1.5]{L}, there is a bijection between the set of $G$-orbits
in $\wfkF\times\wfkF$ and the set $\cnr$ by sending $(L,L')$ to
$A=(a_{i,j})_{i,j\in\bbz}$, where
\begin{equation}\label{aij'} a_{i,j}=\dim_\bbf\frac{L_i\cap
L_j'}{L_{i-1}\cap L_j'+L_i\cap L_{j-1}'}\;\;\text{for $
i,j\in\bbz$},
\end{equation}

\begin{equation}\label{ba}\Theta_\vartriangle(n,r):=\{A=(a_{i,j})_{i,j\in\bbz}\in
M_{\vartriangle,n}(\bbn)\mid  \sum_{1\leq i\leq
n\atop{j\in\bbz}}a_{i,j}=\sum_{1\leq j\leq
n\atop{i\in\bbz}}a_{i,j}=r \}\end{equation} and
$M_{\vartriangle,n}(\bbn)$ is the set of all $\bbz\times\bbz$
 matrices $A=(a_{i,j})_{i,j\in\bbz}$ with $a_{i,j}\in\bbn$ such that
\begin{enumerate}
\item[(a)] $a_{i,j}=a_{i+n,j+n} \;\mbox{for}\; i, j \in \bbz\;; $

\item[(b)]$\mbox{for every}\; i\in \bbz,\;\mbox{the
set}\;\{j\in\bbz\mid a_{i,j}\neq 0\}\;\mbox{is finite}\;. $
\end{enumerate}

Let ${\mathcal O}_A$ denote the orbit in $\wfkF\times\wfkF$
corresponding to $A$. If $(L,L')\in{\mathcal O}_A$, then
$\row(A)=\udim L$ and $\col(A)=\udim L'$, where
$$\row(A)=(\sum_{j\in\bbz}
a_{i,j})_{i\in\bbz},\;\;\col(A)=(\sum_{i\in\bbz}
a_{i,j})_{j\in\bbz}.$$

Assume now that $\bbf=\bbf_q$ is a finite field of $q$ elements
and write $\wfkF(q)$ for $\wfkF$. Suppose that $A, A', A''\in
\cnr$. For any fixed $(L,L'')\in{\mathcal O}_{A''}$, let
$$c_{A,A',A'';q}=|\{L'\in\wfkF(q)\mid (L,L')\in{\mathcal O}_A,
(L',L'')\in{\mathcal O}_{A'}\}|.$$
Clearly, $c_{A,A',A'';q}$ is independent of the choice of $(L,L'')$, and a necessary condition
for $c_{A,A',A'';q}\neq 0$ is that
\begin{equation}\label{ee}
\col(A)=\row(A'),~ \row(A)=\row(A'') ~\text{and}~\col(A')=\col(A'').
\end{equation}

By \cite{L}, there is a polynomial $p_{A,A',A''}\in \zv$ in $v^2$
such that for each finite field $\bbf_q$ with $q$ elements,
$c_{A,A',A'';q}=p_{A,A',A''}|_{v^2=q}$.

\begin{defn} [{\cite{L}},\cite{GV}] \label{def-q-Schur}
The affine quantum Schur algebra $\whS_v(n,r)$ is the free $\zv$-module
with basis $\{e_A\mid A\in\cnr\}$, and
multiplication defined by $$e_{A}\cdot e_{A'}=\begin{cases}\sum_{A''\in\cnr}p_{A,A',A''}e_{A''},~~\text{if}~ \col(A)=\row(A'),\\
              0,~~\text{otherwise}.
\end{cases} $$
\end{defn}

As in the finite case, for each $\lb\in\lr$, define
\diag$(\lb)=(\delta_{i,j}\lb_i)_{i,j\in\bbz}\in
\Theta_\vartriangle(n,r)$, and ${\frak l}_\lb=e_{\text{\diag}(\lb)}$. It is easy to see that for each
$A\in\cnr$,
\begin{equation}\label{idempotents}
{\frak l}_{\lb}e_A=\left\{\begin{array}{ll}
                      e_A,\;&\text{if $\lb=\row(A)$};\\
                     0,
                     &\text{otherwise}\end{array}\right.\;\;\text{and}\;\;
e_A{\frak l}_{\lb}=\left\{\begin{array}{ll}
                      e_A,\;&\text{if $\lb=\col(A)$};\\
                     0, &\text{otherwise.}\end{array}\right.
\end{equation}
 Thus, $\sum_{\lb\in\lr}\frak l_{\lb}$ is the identity of $\whS_v(n,r)$.

For each ring $R$ which is a $\zv$-module, we set
$$\whS(n,r)_R:=\whS_v(n,r)\otimes_{\zv} R.$$
In particular, for the ring of integers $\bbz$ and any field $\mathbb F$, we have
$$\whS(n,r)=\whS(n,r)_\bbz\;\;\text{ and }\;\; \whS(n,r)_{\mathbb F},$$
 where $\bbz$ and $\mathbb F$ are viewed as a $\bbz[v, v^{-1}]$-module by specializing $v$
to $1$.


Now we introduce the algebraic definition of affine Schur algebras
given in \cite{Ya1}. Let $\fkS_r$ denote the symmetric group on
$r$ letters and let $\fsr=\fkS_r\ltimes \bbz^r$ denote the
extended affine Weyl group of type ${\widehat A}_{r-1}$. Define
$$I(n,r)=\{\underline{i}=(i_1,\ldots,i_r)~|~1\leq i_t\leq n, \;\text{for}\; 1\leq t\leq r\},$$
$$I(\bbz,r)=\{\underline{i}=(i_1,\ldots,i_r)~|~i_t\in\bbz, \;\text{for}\; 1\leq t\leq r\}.$$

The group $\fkS_r$ acts on $I(n,r)$ by place permutation while $\fsr$ acts on
$I(\bbz,r)$ on the right with $\fkS_r$ acting by place permutation
and $\bbz^r$ acting by shifting, i.e.,
$$\underline{i}(\sigma,\varepsilon)=\underline{i}+n\varepsilon,$$
for $\underline{i}\in I(\bbz,r),$ $\sigma\in\fkS_r$ and $\varepsilon\in \bbz^r$. Then $\fsr$
acts diagonally on $I(\bbz,r)\times I(\bbz,r)$. For
$(\underline{i}, \underline{j})\in I(\bbz,r)\times I(\bbz,r)$ and
$(\underline{k}, \underline{l})\in I(\bbz,r)\times I(\bbz,r)$, we
identify $\xi_{\underline{i}, \underline{j}}$ and
$\xi_{\underline{k}, \underline{l}}$ if and only if
$(\underline{i}, \underline{j})\sim _\fsr(\underline{k},
\underline{l})$, i.e., $(\underline{i}, \underline{j})$ and
$(\underline{k}, \underline{l})$ are in the same orbit.

\begin{defn} [{\cite{Ya1}}] \label{def-q-Schur2}
The affine Schur algebra $\whS'(n,r)$ is defined to be the $\bbz$-algebra
with basis $\{\xi_{\underline{i}, \underline{j}}~|~\underline{i},
\underline{j} \in I(\bbz, r)/\sim_\fsr\}$ and multiplication given by the
following rule:

$$\xi_{\underline{i},
\underline{j}}\xi_{\underline{k},
\underline{l}}=\sum_{(\underline{p}, \underline{q})\in
I(\bbz,r)\times I(\bbz,r)/\fsr }C(\underline{i},
\underline{j},\underline{k}, \underline{l},\underline{p},
\underline{q})\xi_{\underline{p}, \underline{q}},$$ where
$C(\underline{i}, \underline{j},\underline{k},
\underline{l},\underline{p}, \underline{q})=|\{\underline{s}\in
I(\bbz, r)\mid (\underline{i}, \underline{j})\sim
_\fsr(\underline{p}, \underline{s}),~(\underline{s},
\underline{q})\sim _\fsr(\underline{k}, \underline{l})\}|.$

\end{defn}

\begin{rem}\label{r1} There is a $\bbz$-algebra isomorphism
$$\varphi: \;\; \whS'(n,r)\longrightarrow\whS(n,r)=\whS(n,r)_\bbz,\;\; \xi_{\underline{i}, \underline{j}}\longmapsto e_A=e_A\otimes 1,$$
 where $\ui, \uj\in I(\bbz, r)$ and $A=(a_{x,y})_{x,y\in\bbz}\in\cnr$ is defined by setting
$$a_{x,y}=\sharp\{s~|~i_s=x, j_s=y, 1\leq s\leq r\}.$$
\end{rem}

In the following we will simply identify $\whS'(n,r)$ with $\whS(n,r)$. We now
give some useful multiplication formulas from \cite{Ya1,Ya}.

\begin{prop}[\cite{Ya1}]\label{mu2} Let $\ui, \uj, \underline{k}, \ul\in I(\bbz,
r)$. We have the following equalities in $\sn$.

\begin{itemize}
\item[(1)] $\xi_{\ui, \uj }\xi_{\underline{k}, \ul}=0$ unless $\uj\sim_\fsr
\underline{k}$.

\item[(2)] $\xi_{\ui, \ui }\xi_{\underline{i}, \uj}=\xi_{\underline{i},
\uj}=\xi_{\ui, \uj }\xi_{\underline{j}, \uj}$.

\item[(3)] $\sum_{\ui\in I(n,r)/\fkS_r}\xi_{\ui, \ui }$ is a decomposition
of unity into orthogonal idempotents.
\end{itemize}
\end{prop}

\begin{prop}[\cite{Ya1}, \cite{Ya}]\label{mu4}
\begin{itemize}
\item[(1)]
$$\aligned
\xi_{\underline{i},\underline{j}}\xi_{\underline{j},
\underline{l}}=\sum_{\delta\in\fjl\setminus \fj/ \fij
}\biggl[\fil:\fijl\biggr]\xi_{\underline{i},
\underline{l}\delta}\\
=\sum_{\delta\in\fij\setminus \fj/ \fjl
}\biggl[\ff:\fg\biggr]\xi_{\underline{i}\delta, \underline{l}},\\
\endaligned$$
 where $\ui, \uj, \ul$ are in $I(\bbz,r)$,
$\fii$ is the stabilizer subgroup of $\ui$ in $\fsr$ and $\fij$ is
the stabilizer of $\ui$ and $\uj$ in $\fsr$, i.e. $\fij=\fii\cap
\fj$, etc, $\fjl\setminus \fj/ \fij$ denotes a representative set
of double cosets.

\item[(2)]
$$\aligned
\xij\xjl=\sum_{\delta\in\ffjl\setminus \ffj/ \fflj
}\biggl[\ww:\vv\biggr]\xil\\
 =\sum_{\delta\in\fflj\setminus \ffj/ \ffjl
}\biggl[\wv:\vw\biggr]\yil\\
\endaligned$$
 where $\ui, \uj, \ul$ are in $I(n,r)$, $\varepsilon,
 \varepsilon'\in \bbz^r$,
$\fiii$ is the stabilizer subgroup of $\ui$ in $\fkS_r$ and
$\fijff$ is the stabilizer of $\ui$ and $\uj$ in $\fkS_r$, i.e.
$\fijff=\fiii\cap \fjf$, etc, $\fjlf\setminus \fjf/ \fijf$ denotes
a representative set of double cosets.
\end{itemize}
\end{prop}

\section{canonical basis and two-sided cells of $\whS_v(n,r)$}

In this section we introduce the canonical basis of $\whS_v(n,r)$ given in \cite{L}
and its two-sided cells given in \cite{M}. We finally give an affine
cell chain of $\whS_v(n,r)$ constructed in \cite{C}.

For $A=(a_{i,j})\in \cnr$, let $$[A]=v^{-d_A}e_A,
\;\;\;\;\text{where}\;\; d_A=\sum_{ 1\leq i\leq n\atop{ i\geq k,
j<l,} } a_{i,j}a_{k,l}.$$

Suppose that $A, A_1\in \cnr$ satisfy
$$\row(A)=\row(A_1)=\lb\;\text{ and }\;\col(A)=\col(A_1)=\lb'.$$
 Fix a cyclic flag $L\in \frak F_\lb$, define $A_1\leq A$ if $X_{A_1}^L\subset \bar{X}_{A}^L $, where
$$X_A^L=\{L'\in{\frak F}_{\lb'}\;|\; (L, L')\in {\mathcal O}_A\}.$$

Lusztig has defined in \cite[Sect.~4]{L} the canonical basis
$\{\{A\}| A\in\cnr\}$ of $\whS_v(n,r)$ by
\begin{equation}\label{basis}\{A\}=\sum_{A_1: A_1\leq A}\Pi_{A_1, A}[A_1],
\end{equation}
 where $\Pi_{A, A}=1$ and $\Pi_{A_1, A}\in v^{-1}\bbz[v^{-1}]$ if $A_1<A$.

Recall that for $\lb\in \Lambda^+(n,r)$, ${\flb}=e_{\diag(\lb)}$. We
use $\{\flb\}$ to denote the canonical basis $\{\diag(\lb)\}$.

\begin{lem}\label{10}

For $\lb\in \Lambda^+(n,r)$, we have $\{\flb\}=\flb$ in $\whS_v(n,r)$.
\end{lem}

\begin{pf}

Fix a cyclic flag $L\in\frak F_\lb$. Since ${\flb}=e_{\diag(\lb)}$,
it follows that
$$X_{\diag(\lb)}^L=\{L'\in{\frak F}_{\lb}\;|\; (L, L')\in
{\mathcal O}_{\diag(\lb)}\}=\{L\}.$$
 Applying \eqref{basis} implies that the only term in $\{\flb\}$ with nonzero
 coefficient is $\flb$, that is, $\{\flb\}=\flb$.

\end{pf}

\begin{defn}[\cite{M}]

 For $A, A'\in\cnr$, let
$$\{A\}\{A'\}=\sum_{A''}\nu_{A,A'}^{A''}\{A''\},$$
where $\nu_{A,A'}^{A''}\in\zv$.

 We say that
$A\preceq_{LR}A'$ if there is a sequence $A'=A_1, A_2,\ldots,
A_m=A\in\cnr$ and a sequence $B_1, \ldots, B_{m-1}\in\cnr$ such
that $\nu_{B_s,A_s}^{A_{s+1}}\neq 0$ or
$\nu_{A_s,B_s}^{A_{s+1}}\neq 0$ for $1\leq s\leq m-1$. We write
$A\sim_{LR}A'$ if $A\preceq_{LR}A'$ and $A'\preceq_{LR}A$. The
equivalence classes of $\sim_{LR}$ are called two-sided cells of
$\whS_v(n,r)$.
\end{defn}

Recall that
$\Lambda(n,r):=\{\lambda=(\lambda_1,\ldots,\lambda_n)~|$$~\lambda_i\in\bbn,~$$
\sum_{1\leq i\leq n}\lambda_i=r\}$ and
$$\Lambda^+(n,r):=\{\lambda\in \Lambda(n,r)\mid \lambda_1\geq
\ldots \geq \lambda_n\}.$$

\begin{prop}[\cite{M}]\label{M}

Let $A\in\cnr$. An antidiagonal path in $A$ is an infinite strip
of entries $(a_{i_k,j_k }: k \in\bbz)$ such that $(i_k, j_k)$ is
obtained from $(i_{k-1}, j_{k-1})$ by subtracting $1$ from the
first entry or adding $1$ to the second, with the latter being the
case for all but finitely many $k$. Thus, visually if you draw the
matrix with rows increasing from top to bottom and columns from
left to right (as we will do), then path starts and ends with
infinite vertical strips, and takes finitely many right or
vertical turns.

Let $d_j$ be the maximal size of the sum of entries in the union
of $j$ antidiagonal paths. Then, we define a
map\begin{equation}\label{map rho}\aligned\rho:
\cnr&\longrightarrow \Lambda^+(n,r),\\
A&\longmapsto (d_1, d_2-d_1,\ldots, d_n-d_{n-1}).
\endaligned\end{equation}
Define
$${\frak c}_\lb=\{\{A\}\;|A\in \rho^{-1}(\lb)\}.$$
 Then $\{{\frak c}_\lb\;|\;\lb\in\Lambda^+(n,r)\} $ are the two-sided cells of $\whS_v(n,r)$.

\end{prop}

\begin{lem}\label{00}

For $\lb\in \Lambda^+(n,r)$, we have $\{\flb\}\in \frak c_\lb$.

\end{lem}

\begin{pf}

This can be directly checked by using Proposition \ref{M}.

\end{pf}

Let ``$\geq$'' be a partial ordering on $\Lambda^+(n,r)$ defined by
$\lambda \geq \mu$ whenever $\sum_{j=1}^i \lambda_j\geq
\sum_{j=1}^i \mu_j$, for $1\leq i\leq n$. It is obvious that
$(r,0,\ldots,0)$ is the maximal element with respect to this ordering.
We then fix a total ordering
$\lb^{(1)}>\lb^{(2)}>\ldots>\lb^{(t)}$ on $\Lambda^+(n,r)$ which is a refinement of the above ordering,
where $t$ equals the number of elements in $\Lambda^+(n,r)$.

For each $1\leq i\leq t$, define
$${\frak c}_i=\{\{A\}\mid A\in \rho^{-1}(\lb^{(i)})\}.$$
 Then $\frak c_1, \frak c_2,\ldots, \frak c_t$ are the two sided
cells in $\whS_v(n,r)$ such that $\frak c_i$$\preccurlyeq_{LR} \frak c_j$
implies $i\leq j$.

For $1\leq i\leq t$, let $C_i'$ be the $\bbz[v,v^{-1}]$-submodule of $\whS_v(n,r)$
generated by all $\{e_A\}$ with $A\in \frak c_i$. Let
$C_i=\oplus_{l=1}^iC_l'$. Then by \cite{C}, $C_i$ is an ideal of
$\whS_v(n,r)$ for each $1\leq i\leq t$ and the chain of ideals
\begin{equation}\label{cui}
0\subseteq C_1\subseteq \ldots \subseteq C_i\subseteq\ldots \subseteq C_t=\whS_v(n,r)
\end{equation}
forms an affine cell chain of $\whS_v(n,r)$. This implies particularly that $\whS_v(n,r)$
is an affine cellular algebra over $\bbz[v,v^{-1}]$.

Now applying Lemma \ref{ring} to the evaluation map $\bbz[v,v^{-1}]\rightarrow \bbz$ taking $v^{\pm}\mapsto 1$,
we conclude that $\whS(n,r)$ is an affine cellular $\bbz$-algebra which admits an affine cell chain
induced from \eqref{cui}. Thus, in the following sections we simply view the $C_i$ as ideals of $\whS(n,r)$ and the
chain \eqref{cui} as an affine cell chain of $\whS(n,r)$.

\section{Affine quasi-heredity of $\whS(n,r)$ in case $n> r$}

In this section we prove that $\whS(n,r)$ is affine
quasi-hereditary under the assumption $n>r$. Keep all the notations in the previous
sections.

Since $\whS(n,r)$ is a kind of BLN algebra when $n>r$, it follows from \cite{N}
(see also \cite{C1}) that affine cell ideals of $\whS(n,r)$ are idempotent
ideals.

\begin{prop}[\cite{N}, \cite{C1}]\label{id ideal}
Suppose that $n>r$. Then for each $1\leq i\leq t$,
$C_i'=C_i/C_{i-1}$ is an idempotent ideal of $\whS(n,r)/C_{i-1}$.
\end{prop}

Fix a total ordering
$\lb^{(1)}>\lb^{(2)}>\ldots>\lb^{(t)}$ on $\Lambda^+(n,r)$ as in Section 4,
where $t$ equals the number of elements in $\Lambda^+(n,r)$.

For each $1\leq i\leq t$, set
$$J_{i}=\whS(n,r)\big(\sum_{\mu\in
\Lambda^+(n,r),\;\mu\geq\lb^{(i)}}{\frak l}_\mu\big)\whS(n,r).$$

\begin{lem}\label{11}
 Assume $n>r$. Then $J_i= C_i$ for each $1\leq i\leq t$. Moreover,
 $\bar{J}_i=J_i/J_{i-1}$ is projective as an
$\bar{S}$-module, where $\bar{S}=\whS(n,r)/J_{i-1}$.

\end{lem}

\begin{pf} First we remark that $\lb^{(1)}=(r,0,\ldots,0)\in\Lambda^+(n,r)$. By Lemmas
\ref{10} and \ref{00},
$${\fl}=\{{\fl}\}\in C_1'.$$
 Since $C_1$ is an idempotent ideal of $\whS(n,r)$ by Proposition \ref{id ideal}, it follows from Lemma \ref{usu} that
$$C_1=\whS(n,r) {\fl}\whS(n,r)=J_1,$$
and $J_1$ is projective as a left $\whS(n,r)$-module.

For each $1\leq i\leq t$, the above analysis is valid for $\bs=\whS(n,r)/J_{i-1}$.
Then we get that
\begin{equation}\label{ide}
C_i'=C_i/C_{i-1}=\bs{\bfli}\bs,
\end{equation}
 and $C_i/C_{i-1}$ is projective as a left
$\bs$-module.

Now we prove that $C_i=J_i$ for $2\leq i\leq t$ by induction. By
Lemmas \ref{10} and \ref{00},
$$\sum_{1\leq l\leq i}{\frak l}_{\lb^{(l)}}=\sum_{1\leq l\leq i}\{{\frak l}_{\lb^{(l)}}\}\in C_i.$$
Then there is a canonical inclusion from $J_i$ to $C_i$
for $1\leq i\leq t$. This together with the induction hypothesis induces an inclusion of $\bbz$-modules
$$\theta:\;\; J_i/J_{i-1}\hookrightarrow C_i/J_{i-1}=C_i/C_{i-1}$$
 taking $\bfli \mapsto \bfli$. It is easy to check that
$$J_i/J_{i-1}\cong \bs\bfli\bs.$$  On the other hand,
$C_i/C_{i-1}= \bs\bfli\bs$ by \eqref{ide}. Then $\theta$ is an
isomorphism. This proves that $J_i= C_i$.
\end{pf}

Now we introduce the representation ring of the general linear
group. It has been shown in \cite{C} that the affine cell ideals
of $\whS(n,r)$ can be considered as generalized matrix algebras over
these representation rings.

Let $R(\GL_m(\bbc))$ denote the representation ring of $\GL_m(\bbc)$,
where $m\in\bbn$. It is shown in \cite[Exercise 23.36]{FH} that
$$R(\GL_{m}(\bbc))\cong \bbz[x_1, x_2,\ldots, x_{m}, x_{m}^{-1}].$$

For $\lb=(\lb_1,\ldots,\lb_n)\in \Lambda^+(n,r)$, denote by
$R(\prod_{l=1}^n\GL_{\lb(l)}(\bbc))$ the representation ring of
$\prod_{l=1}^n\GL_{\lb(l)}(\bbc)$, where $\lb(l)=\lb_l-\lb_{l+1}$
with $\lb_{n+1}=0$. Define
\begin{equation}\label{bi}B_\lb=R\big(\prod_{l=1}^n\GL_{\lb(l)}(\bbc)\big).\end{equation}
 Then
$$B_\lb=\bbz[x_1, x_2, \ldots, x_{\lb_1}, x_{m(1)}^{-1}, x_{m(2)}^{-1},\ldots,x_{m(n)}^{-1}] ,$$
where $m(i)=\lb_1-\lb_{i+1}$ and $\lb_{n+1}=0$. For example, if $\lb=(4, 2, 1)$, then
$$B_\lb=\bbz[x_1, x_2, x_2^{-1}, x_3, x_3^{-1}, x_4, x_4^{-1}].$$

\begin{thm}\label{main}
In case $n>r$, $\whS(n,r)$ is an affine hereditary $\bbz$-algebra.
\end{thm}

\begin{pf}  By \cite[Th.~5]{y}, $\whS(n,r)$ is Noetherian.
Recall that $\lb^{(1)}>\lb^{(2)}>\ldots>\lb^{(t)}$ is a total
ordering on the elements of $\Lambda^+(n,r)$ and
$$J_{i}=\whS(n,r)\big(\sum_{\mu\in
\Lambda^+(n,r),\;\mu\geq\lb^{(i)}}{\frak l}_\mu\big)\whS(n,r),$$
 for $1\leq i\leq t$. Then we obtain the following chain of ideals of $\whS(n,r)$:
$$0=J_0\subseteq J_1\subseteq \ldots \subseteq J_t=\whS(n,r).$$
 It suffices to prove that $J_i/J_{i-1}$ is an affine-heredity
ideal of $\bs=\whS(n,r)/J_{i-1}$ for $1\leq i\leq t$.
 By Lemma \ref{11}, $J_i/J_{i-1}=\bs\bfli$$\bs$ is an idempotent  ideal of
 $\bar{S}$ and $\bs\bfli\bs$ is projective over $\bar{S}$.
 Thus, condition $(1)$ in Definition \ref{quasi} is fulfilled.

By Definition \ref{def-aff} and Lemma \ref{11}, for each $1\leq i\leq
t$, there are $\bs$-$B_{\lb^{(i)}}$-bimodule $\Delta$ and
$B_{\lb^{(i)}}$-$\bs$-bimodule $\Delta'$ with the following
$\bs$-bimodule isomorphism
$$\alpha:J_i/J_{i-1}\longrightarrow \Delta\otimes_{B_{\lb^{(i)}}}\Delta',$$
where $B_{\lb^{(i)}}$ is defined as in \eqref{bi}.

For notational simplicity, we write $B_i$ for $B_{\lb^{(i)}}$. We remark that $\Delta$ and $\Delta'$ are
free of the same finite rank over $B_i$ when viewed as right and
left modules, respectively.

Since $ J_i/J_{i-1}=\bs \bfli\bs$ is generated by the idempotent
$\bfli$, it follows from Lemma \ref{usu} that
$$\add(_{\bs}\bs \bfli)=\add(_{\bs}\Delta).$$
This implies that $\Delta$ is projective as a left $\bs$-module. Since $\Delta'$ is a
free $B_i$-module of finite rank, we get the following
decomposition of $\bs\bfli\bs$ as a direct sum of projective
$\bs$-modules:
$$ \bs\bfli\bs=J_i/J_{i-1}\cong \underbrace{\Delta \oplus\ldots\oplus \Delta}_m,$$
where $m$ equals the rank of $\Delta'$ over $B_i$.

To prove that $ \bs\bfli\bs$ satisfies condition (2) of Definition
\ref{quasi}, we only need to show that $\End_{\bs}(\Delta)$ is an
affine algebra. Indeed, by the definition of $B_i$, it is easy to see that
$\rad(B_i)=0$. Then, by Lemma \ref{usu}, $\End_{\bs}(\Delta)\cong B_i$ is an affine algebra.

Since $\Delta$ is free of finite rank over $B_i$, condition (3) in Definition \ref{quasi} is satisfied. This
finishes the proof.
\end{pf}

\section{Affine quasi-heredity of $\whS(n,r)$ in case $n\leq r$}

In this section we prove that $\whS(n,r)$ is affine quasi-hereditary when $n\leq r$.
We first give some properties of $\whS(n,r)$ which will be needed in the proof.
Note that these properties are valid for $\whS(n,r)$ with $n,r$ arbitrary.

Recall that there is an action of the symmetric group $\fkS_n$ on
$\Lambda(n,r)$ given by
\begin{equation}\label{act}w\cdot(\lambda_1,\ldots,\lambda_n)=(\lambda_{w(1)},\ldots,\lambda_{w(n)}),\end{equation}
where $w\in \fkS_n$.

\begin{lem}\label{22}
Let $\lambda, \mu\in \Lambda(n, r)$. If there is some $w\in
\fkS_n$ such that $w\cdot \lambda=\mu$, then
$$\whS(n,r) {\frak l}_\lambda \whS(n,r)=\whS(n,r){\frak l}_\mu\whS(n,r).$$
\end{lem}

\begin{pf} By \cite[Lem.~2]{Ya}, $\whS(n,r){\frak l}_\lambda \cong \whS(n,r){\frak l}_\mu$.
The isomorphism is induced by two elements $X, Y\in\whS(n,r)$ such
that
$$X\in {\frak l}_\lambda \whS(n,r) {\frak l}_\mu ,\; Y\in {\frak l}_\mu \whS(n,r)
{\frak l}_\lambda ,\;  XY={\frak l}_\lambda, \;\text{ and }\; YX= {\frak l}_\mu.$$
 Since $X{\frak l}_\mu=X$ and $Y{\frak l}_\lambda=Y$, we get that
$$X{\frak l}_\mu Y=XY={\frak l}_\lambda \;\text{ and }\; Y{\frak l}_\lambda X=YX={\frak l}_\mu.$$
This implies that ${\frak l}_\lambda$ is contained in the ideal
generated by ${\frak l}_\mu$, and vice versa. Therefore,
$$\whS(n,r) {\frak l}_\lambda \whS(n,r)=\whS(n,r) {\frak l}_\mu \whS(n,r).$$
\end{pf}

\begin{lem}Write $\whS=\whS(n,r)$. Then
$$\whS=\whS\big(\sum_{\lambda\in
\Lambda(n,r)} {\frak l}_\lambda\big) \whS=\whS\big(\sum_{\lambda\in
\Lambda^+(n,r)}{\frak l}_\lb\big)\whS.$$

\end{lem}

\begin{pf} Since
$$\sum_{\lambda\in \Lambda(n,r)}{\frak l}_\lambda={\rm id},$$
we get that
$$\whS=\whS(\sum_{\lambda\in\Lambda(n,r)}{\frak l}_\lb)\whS.$$

By \eqref{act}, each orbit of $\Lambda(n,r)$ under the action of
$\fkS_n$ has exactly one element in $\Lambda^+(n,r)$. It follows from
Lemma \ref{22} that
$$\whS=\whS\big(\sum_{\lambda\in
\Lambda(n,r)} {\frak l}_\lambda\big) \whS=\whS\big(\sum_{\lambda\in
\Lambda^+(n,r)}{\frak l}_\lb\big)\whS,$$
 as desired.
\end{pf}

Now we are ready to prove the main result in this section. Choose $N>r$. Then by Theorem \ref{main},
$\whS(N,r)$ is affine quasi-hereditary with an affine-heredity chain
\begin{equation}\label{41}0=J_0\subseteq J_1\subseteq\cdots J_i\subseteq\cdots\subseteq
J_t=\whS(N,r),\end{equation} where $t$ is the number of elements in
$\Lambda^+(N,r)$ and
$$J_i=\whS(N,r)\big(\sum _{\mu\in \Lambda^+(N,r), \;\mu\geq\lb^{(i)}}{\frak l}_\mu\big)\whS(N,r).$$

Since $n\leq r<N$, there is an inclusion map from $\Lambda(n,r)$ to $\Lambda(N,r)$
$$\aligned \Lambda(n,r)&\longrightarrow \Lambda(N,r)\\
(\lb_1,\ldots,\lb_n)&\longmapsto
(\lb_1,\ldots,\lb_n,0,\ldots,0).\endaligned$$
 Then we can view $\Lambda(n,r)$ as a subset of $\Lambda(N,r)$. Furthermore, let
$$e=\sum_{\lb\in \Lambda(n,r)}\frak l_\lb.$$
Then $e^2=e\in \whS(N,r)$ and it is easily checked that
\begin{equation}\label{eAe}
e\whS(N,r) e\cong \whS(n,r).
\end{equation}

\begin{lem}\label{23}
Let $\{J_i\;|\;1\leq i\leq t\}$ be the ideals of $\whS(N,r)$ in the
affine quasi-heredity chain as given in \eqref{41}. Then there is
some $1\leq l \leq t$ such that $J_l=\whS(N,r) e\whS(N,r)$.
\end{lem}

\begin{pf} Set
$$x=\sum_{\lb\in \Lambda^+(n,r)}\frak l_\lb.$$
 Then there is some $1\leq l \leq m$ such that $J_l=\whS(N,r) x\whS(N,r)$.
By Lemma \ref{22},
$$\whS(N,r) e\whS(N,r)=\whS(N,r) x\whS(N,r)=J_l.$$

\end{pf}

\begin{thm}\label{24} The algebra $e\whS(N,r) e$ is affine quasi-hereditary. In other words,
$\whS(n,r)\cong e\whS(N,r)e$ is affine quasi-hereditary when $n\leq r$.
\end{thm}

\begin{pf} Take the affine-heredity chain of $\whS(N,r)$ as in \eqref{41}
$$0=J_0\subseteq J_1\subseteq\cdots J_l\subseteq\cdots\subseteq J_t=\whS(N,r)$$
such that $J_l=\whS(N,r) e\whS(N,r)$, where $e=\sum_{\lb\in \Lambda(n,r)}\frak l_\lb$.

We want to show that $$0=eJ_0e\subseteq eJ_1e\subseteq\cdots
eJ_ie\subseteq\cdots\subseteq eJ_le=e\whS(N,r) e$$ is an affine-heredity
chain of $e\whS(N,r) e$, i.e., $\whS(n,r)=e\whS(N,r) e$ is affine quasi-hereditary.

First we prove that $eJ_le=e\whS(N,r) e$. To simplify the notation, we
use $\whS$ to denote $\whS(N,r)$ in what follows. Since $J_l=\whS e\whS$,
we deduce that
$$eJ_le=e\whS e\whS e=e(\whS e)^2=e\whS e.$$

Now we only need to prove that $eJ_1e$ is an affine-hereditary
ideal of $e\whS e$, the other requirements follow easily by
induction.

We first check that $eJ_1e=(eJ_1e)^2$. It is obvious that
$eJ_1e\supseteq(eJ_1e)^2$. Note that $J_1=\whS{\frak l}_{\lb^{(1)}} \whS$,
where $\lb^{(1)}=(r,0,\ldots,0)\in \Lambda^+(N,r)$. Since
$$\aligned
(eJ_1e)^2&=eJ_1eJ_1e\\
&=e\whS{\frak l}_{\lb^{(1)}} \whS e\whS{\frak l}_{\lb^{(1)}}
\whS e\;\;\;\;\;\;\text{by}\; J_1=\whS{\frak l}_{\lb^{(1)}} \whS\\
&\supseteq e\whS{\frak l}_{\lb^{(1)}} \whS e\;\;\;\;\;\;\;\;\;\;\;\;\;\;\;\;\;\text{by}\; {\frak
l}_{\lb^{(1)}} \in \whS e\whS\\
&=eJ_1e,\endaligned$$ we get that $eJ_1e=(eJ_1e)^2$. This proves
that $eJ_1e$ satisfies condition $(1)$ of Definition \ref{quasi}.

Now we want to show that $eJ_1e$ can be decomposed into a direct sum
of some projective $e\whS e$-module $P$ and that the endomorphism algebra
of $P$ is an affine algebra.

Since $J_1$ is an affine cell ideal of $\whS$, there are
$\whS$-$B$-bimodule $\Delta$ and $B$-$\whS$-bimodule $\Delta'$ with the
following $\whS$-$\whS$-bimodule isomorphism
$$J_1\longrightarrow\Delta\otimes_{B}\Delta',$$  where
$B=B_{\lb^{(1)}}$, $\Delta$ and $\Delta'$ are free over $B$ of
finite rank. By \cite[Prop.~6]{Ya},
$$B=B_{\lb^{(1)}}\cong \bbz[x_1, \ldots, x_r, x_r^{-1}].$$ Since
$J_1$ is idempotent and contains a nonzero idempotent ${\frak l}_{\lb^{(1)}} $, we get
by Lemma \ref{usu} that
$$\add(\whS{\frak l}_{\lb^{(1)}})=\add(\Delta).$$
Then $\Delta$ is projective as a left $\whS$-module
and $J_1$ can be decomposed into a direct sum of projective
$\whS$-modules
$$J_1\cong\underbrace{\Delta\oplus\ldots\oplus \Delta}_{s},$$
where $s$ is the rank of $\Delta'$ over $B$.

Further, by \cite[Lem.~3.3]{y1},
\begin{equation}\label{717}
eJ_1e\longrightarrow(e\Delta)\otimes_{B}(\Delta'e)
\end{equation}
is an $e\whS e$-bimodule isomorphism which makes $eJ_1e$ an affine cell
ideal of $e\whS e$. Since $eJ_1e$ is idempotent and contains a
nonzero idempotent ${\frak l}_{\lb^{(1)}}$, it follows from Lemma \ref{usu} that
\begin{equation}\label{716}
\add\biggl(_{e\whS e}(e\whS e){\frak l}_{\lb^{(1)}}\biggr)=\add\biggl(_{e\whS e}e\Delta\biggr).
\end{equation}
 Hence, $e\Delta$ is a left projective $e\whS e$-module. Since $\Delta' e$
 is free of finite rank  over $B$, we obtain an $e\whS e$-module decomposition
$$eJ_1e=\underbrace{e\Delta \oplus\ldots\oplus e \Delta }_{a},$$
where $a$ is the rank of $\Delta' e$ over $B$.

Since $\rad(B)=0$, we get that $\End_{e\whS e}(e\Delta )\cong B$ by
Lemma \ref{usu}. Thus, $\End_{e\whS e}(e\Delta )$ is an affine algebra.
Consequently, $eJ_1e$ satisfies condition (2) of Definition \ref{quasi}.

Since $\Delta $ is free of finite rank over $B$ as a right module, it is easy to check that $e\Delta$ is projective and
finitely generated over $B$. Thus, condition (3) of
Definition \ref{quasi} is satisfied. Hence, $eJ_1e$ is an affine-hereditary ideal
of $e\whS e$. The proof is completed.
\end{pf}

Theorem \ref{24} together with Theorem \ref{main} gives a positive answer to
a conjecture in \cite{K} and \cite{C} which states that $\whS(n,r)$ is affine quasi-hereditary.
Furthermore, the affine-heredity chain of $\whS(n,r)$ satisfies the conditions in
\cite[Th.~4.4]{KX}. i.e., for each $1\leq i\leq t$, $\rad (B_{\lb^{(i)}})=0$, $\bs\bfli\bs$ is
idempotent and contains a nonzero idempotent in $\bs$. We remark that in this case,
the global dimension of $\whS(n,r)$ is finite if and only
if the global dimension of $B_{\lb^{(i)}}$ is finite for every $1\leq i\leq t$. Since $B_{\lb^{(i)}}$ is a
localization of $\bbz[x_1,\ldots,x_{\lb_1}]$, it has finite global dimension. Therefore, we obtain the
following corollary.

\begin{cor}  For all positive integers $n$ and $r$, the global dimension
of $\whS(n,r)$ is finite.
\end{cor}

\section{Application I: A description of subquotient algebras of $\whS(n,r)_\bbq$}

In this section we always work with the $\bbq$-algebra $\whS(n,r)_\bbq$. For
$\lb=(\lb_1,\ldots,\lb_n)\in \Lambda^+(n,r)$, set
\begin{equation}\label{bi1}
B_\lb=(B_\lb)_\bbq=R\big(\prod_{l=1}^n\GL_{\lb(l)}(\bbc)\big)\otimes_\bbz\bbq,
\end{equation}
where $\lb(l)=\lb_l-\lb_{l+1}$ and $\lb_{n+1}=0$. i.e.
$$B_\lb=\bbq[x_1, x_2, \ldots, x_{\lb_1}, x_{m(1)}^{-1}, x_{m(2)}^{-1},\ldots,x_{m(n)}^{-1}] ,$$
where $m(i)=\lb_1-\lb_{i+1}$ and $\lb_{n+1}=0$.

The main aim of this section is to prove that for each $\lb\in \Lambda^+(n,r)$,
$$\bfl\bs\bfl\cong B_\lb,$$
 where $\bs=\whS(n,r)_\bbq/J$ with $J=\whS(n,r)_\bbq(\sum_{\mu\in
\Lambda^+(n,r),\;\mu>\lb}  {\frak l}_\mu)\whS(n,r)_\bbq$.

First we have the following lemma.

\begin{lem}\label{71} For each $\lb\in \Lambda^+(n,r)$, $\bfl\bs\bfl$ is Morita equivalent
to $B_\lb$.
\end{lem}

\begin{pf} It is clear that $\bfl\bs\bfl\cong\End_{\bs}(\bs\bfl)$. Since
$\bs\bfl\bs$ is an affine cell ideal of $\bs$, there are
$\bs$-$B_\lb$-bimodule $\Delta$ and $B_\lb$-$\bs$-bimodule
$\Delta'$ with the following $\bs$-$\bs$-bimodule isomorphism
$$\bs\bfl\bs\longrightarrow \Delta \otimes_{B_\lb}\Delta'.$$
Then $\add_{\bs}(\bs\bfl)=\add_{\bs}(\Delta)$ by Lemma
\ref{usu}. This implies that $\bfl\bs\bfl$ and
$\End_{\bs}(\Delta)$ are Morita equivalent. Since
$\rad{(B_\lb)}=0$, we get that $\End_{\bs}(\Delta)\cong B_\lb$ by
Lemma \ref{usu}. This proves the lemma.

\end{pf}

By the commutativity of $B_\lb$ and Lemma \ref{71}, in order to get the isomorphism
$\bfl\bs\bfl\cong B_\lb$, it suffices to prove that $\bs\bfl\bs$ is commutative. To prove this fact,
we need the following basic results on classical Schur algebra $S(n,r)_\bbq$ (over $\bbq$).

We first recall a refined triangular decomposition and a
presentation of $S(n,r)_\bbq$ given in \cite{DG} and
\cite{DG1}, respectively. For $1\leq i, j\leq n$, let
$E_{i,j}=(a_{k,l})$ be the $n\times n$ elementary matrix in $\bbn$
 defined by
$$a_{k,l}=\begin{cases}1,~~\text{if}~ k=i, l=j,\\
              0,~~\text{otherwise}.
\end{cases}$$ Let
$$e_i=\sum_{\lb\in\Lambda(n,r-1)}e_{E_{i,i+1}+\text{diag}(\lb)},\;
f_i=\sum_{\lb\in\Lambda(n,r-1)}e_{E_{i+1,i}+\text{diag}(\lb)},\;
{\frak l}_\lb=e_{\text{diag}(\lb)}$$
 for $1\leq i<n$ and $\lb\in\Lambda(n,r)$. Suppose that $I$ is the set of finite
sequence $B=(i_1,\ldots,i_m)$ (including the empty sequence
$\emptyset$), with each $i_l$ for $1\leq l\leq m$ contained in the
set $\{1,2,\ldots,n-1\}$. For each sequence $B=(i_1,\ldots,i_m)$,
define
$$e_B=e_{i_1}e_{i_2}\cdots e_{i_m};~~f_B=f_{i_m}\cdots
f_{i_2}f_{i_1}$$ and set $e_\emptyset=f_\emptyset=1$ by convention.

\begin{lem}[{\cite[Prop 2.7]{DG}}]\label{dg1} Let $S(n,r)_\bbq$ be the classical Schur algebra
over $\bbq$. Then $S(n,r)_\bbq$ is spanned by products of the form
$f_Be_\lb e_D$, where $B,D$ are in $I$ and $\lb\in\Lambda(n,r)$.
\end{lem}

\begin{lem}[{\cite[Thm 1.4]{DG1}}]\label{dg2} The $\bbq$-algebra $S(n,r)$ is generated by
$$e_i, f_i, {\frak l}_\lb ~~(1\leq i\leq n-1, \lb\in\Lambda(n,r))$$subject
to the following relations

\begin{enumerate}
\item[(1)] ${\frak l}_\lb {\frak l}_\mu=\delta_{\lb,\mu}{\frak l}_\lb,~
\sum_{\lb\in\Lambda(n,r)}{\frak l}_\lb=1$,

    \item[(2)] $e_i{\frak l}_\lb=
   \begin{cases} {\frak l}_{\lb+\alpha_i-\alpha_{i+1}}e_i,~~\text{if}~\lb_{i+1}\geq 1,\\
              0,~~\lb_{i+1}=0.
\end{cases}$

\item[(3)] ${\frak l}_\lb e_i=
   \begin{cases} e_i{\frak l}_{\lb-\alpha_i+\alpha_{{i+1}}},~~\text{if}~\lb_{i}\geq 1,\\
              0,~~\lb_{i+1}=0.
\end{cases}$

 \item[(4)] $f_i{\frak l}_\lb=
   \begin{cases} {\frak l}_{\lb-\alpha_i+\alpha_{i+1}}f_i,~~\text{if}~\lb_{i}\geq 1,\\
              0,~~\lb_{i}=0.
\end{cases}$

\item[(5)] ${\frak l}_\lb f_i=
   \begin{cases} f_i{\frak l}_{\lb+\alpha_i-\alpha_{i+1}},~~\text{if}~\lb_{i+1}\geq 1,\\
              0,~~\lb_{i}=0.
\end{cases}$

\item[(6)] $e_i
f_j-f_je_i=\delta_{i,j}\sum_{\lb\in\Lambda(n,r)}(\lb_i-\lb_{i+1})\frak
l_\lb,$

\item[(7)] $e_i^2e_{j}-2e_ie_{j}e_i+e_{j}e_i^2=0,$ \;
$(|i-j|=1)$\\
$e_ie_j=e_je_i$ \;(\text{otherwise})

\item[(8)] $f_i^2f_{j}-2f_if_{j}f_i+f_{j}f_i^2=0,$ \;
$(|i-j|=1)$\\
$f_if_j=f_jf_i$ \;(\text{otherwise})

\end{enumerate}
where $\alpha_i=(0,\ldots,\underset{(i)}1,\ldots,0)\in\bbz^n$.

\end{lem}

\begin{thm}\label{iso}  Let $\lb=(\lb_1, \ldots, \lb_n)\in \Lambda^+(n,r)$. Then
$\bfl\bs\bfl\cong B_\lb.$
\end{thm}

\begin{pf} We identify $\frak l_\lb$ with $\xi_{\ui,\ui}$, where
$\ui=(\underbrace{ 1,\ldots,1}_{\lb_1},\ldots,\underbrace{
n,\ldots,n}_{\lb_n})\in I(n,r)$. Then ${\flb}\whS(n,r)_\bbq {\flb}$ is
spanned by the set
$${\frak Y}=\{e_A=\xi_{\ui,\ui \sigma+n\varepsilon }\mid \sigma\in \fkS_r, \varepsilon\in \bbz^r \}.$$
 Let
 $${\frak X}=\{e_A=\xi_{\ui,\ui \sigma}\mid\sigma\in \fkS_r \}$$
 be a subset of $\frak Y$. This implies that
$${\frak X}=\{e_A\mid A=(a_{i,j})_{i,j\in\bbz}\in\cnr, \col(A)=\row(A)=\lb, \sum_{1\leq s,t\leq n}a_{s,t}=r\}.$$
 In this sense, we can view $\frak{X}$ as a subset of ${\flb}
S(n,r){\flb}$. By Lemma \ref{dg1}, for each $e_A\in\frak{X}$, we can
write
$$e_A=\sum_sa_s(f_{B_s}{\frak l}_{\lb^{(s)}} e_{D_s}),$$
where $a_s\in\bbq$, $B_s, D_s\in I$ and $\lb^{(s)}\in\Lambda(n,r)$. Consequently,
$$e_A={\flb} e_A {\flb}=\sum_s{\flb} f_{B_s}{\frak l}_{\lb^{(s)}} e_{D_s}\flb.$$
 If $D_s\neq\emptyset$ and $e_{D_s}\flb\neq 0$ for some $s$, then
$$e_{D_s}{\flb}={\frak l}_{\lb'} e_{D_s},$$ where $\lb'> \lb$ by
Lemma \ref{dg2}. This means that ${\flb} f_{B_s}e_{\lb^{(s)}}
e_{D_s}{\flb}\in J$. Similarly, if $B_s\neq\emptyset$ and ${\flb} f_{B_s}\neq 0$
for some $s$, then ${\flb} f_{B_s}{\frak l}_{\lb^{(s)}} e_{D_s}{\flb}\in J$.
Therefore, for each $A\in \frak X$,
$$\overline{e}_A=\overline{{\flb} e_A {\flb}}=\sum_s\overline{{\flb}
{\frak l}_{\lb^{(s)}}{\flb}}=a\overline{\frak l}_\lb,
$$ where $a$ equals the number of $s$ such that $\lb^{(s)}=\lb$. Hence, $\bfl\bs\bfl$ is $\bbq$-spanned by
\begin{equation}\label{span}{\frak M}=\{\overline{e}_A=\overline{\xi}_{\ui,\ui +n\varepsilon }\mid
\varepsilon\in \bbz^r\}.\end{equation}

Let $\Lambda$ be the subalgebra of ${\frak l}_\lb\whS(n,r)_\bbq{\frak l}_\lb$
which is $\bbq$-spanned by $\{\xi_{\ui,\ui +n\varepsilon }\mid
\varepsilon\in \bbz^r \}$. Applying \cite[Prop.~8]{Ya} gives
$$\Lambda\cong\bbc[x_1, x_2,\ldots,
x_r, x_{1}^{-1}, x_{2}^{-1},\ldots, x_{\alpha}^{-1}],$$ where
$\alpha$ is the number of nonzero entries in $\{\lb_1,
\ldots,\lb_n\}$. Then by \eqref{span}, we get an epimorphism from $\Lambda$ to
$\bfl\bs\bfl$. This proves that $\bfl\bs\bfl$ is a commutative
algebra. Applying Lemma \ref{71} gives ${\bfl}\bs{\bfl}\cong
B_\lb$. This proves the theorem.

\end{pf}

\section{Application II: An identification of parameter sets of simple modules}

In this section we give a parameter set of simple
$\whS(n,r)_\bbc$-modules and further identify this parameter set
with the set given in \cite{DDF} and \cite{F}.

For $\lb=(\lb_1,\ldots,\lb_n)\in \Lambda^+(n,r)$, recall that
$$B_\lb=\bbc[x_1, x_2, \ldots, x_{\lb_1}, x_{m(1)}^{-1}, x_{m(2)}^{-1},\ldots,x_{m(n)}^{-1}] ,$$
where $m(i)=\lb_1-\lb_{i+1}$ and $\lb_{n+1}=0$. Define
$$\Omega_\lb=\{\underline{a}=(a_1, a_2, \ldots, a_{\lambda_1})\in\bbc^{\lambda_1}\;|\;a_{m(i)}\neq 0, \;\text{for}\; 1\leq i\leq n\},$$
where $m(i)=\lb_1-\lb_{i+1}$ and $\lb_{n+1}=0$.

For $\underline{a}=(a_1, a_2, \ldots, a_{\lb_1})\in \Omega_\lb$,
define
$$E_{\underline{a}}=B_\lb/(x_1-a_1, \ldots, x_{\lb_1}-a_{\lb_1} ).$$
Then we have the following lemma.

\begin{lem}\label{80}

For $\lb\in \Lambda^+(n,r)$, the set
$$\{E_{\underline{a}}\;|\; \underline{a}\in\Omega_\lb\}$$
forms a complete set of non-isomorphic simple $B_\lb$-modules.

\end{lem}

\begin{lem}\label{81}

Let $\Omega_{r,n}=\cup_{\lb\in\Lambda^+(n,r) } \Omega_\lb$. Then
$\Omega_{r,n}$ is a parameter set of non-isomorphic simple $\whS(n,r)_\bbc$-modules.
\end{lem}

\begin{pf} Let
$$J_0=0\subseteq J_1\subseteq \cdots\subseteq J_t=\whS(n,r)_\bbc$$
be the cell chain given in Theorem \ref{main} and Theorem
\ref{24} such that $J_i/J_{i-1}\cong \bs\bfli\whS(n,r)_\bbc$, where
$\bs=\whS(n,r)_\bbc/J_{i-1}$.

By \cite[Cor.~3.2]{KX}, simple $\whS(n,r)_\bbc$-modules are
parameterized by simple $J_i/J_{i-1}$-modules for all $1\leq i\leq
t$. Here for an ideal $J_i/J_{i-1}$ of $\bs$, a
$J_i/J_{i-1}$-module $M$ is called simple if $(J_i/J_{i-1})M\neq
0$ and there are no submodules different from $0$ and $M$. Since
$J_i/J_{i-1}\cong \bs\bfli\bs$, we get that the correspondence
$L\mapsto \bfli L$ induces a bijection between the set of
isomorphism classes of simple $J_i/J_{i-1}$-modules and that
 of isomorphism classes of simple $\bfli\bs\bfli$-modules.

By Lemma \ref{iso}, $\bfli\bs\bfli\cong B_{\lb^{(i)}}$. The lemma then follows
from Lemma \ref{80}.
\end{pf}

Now we recall the definition of segments in \cite{DDF} (by specializing $v$ to $1$). A segment
$s$ is by definition a sequence
$$s=(a, a, \ldots, a)\in (\bbc^*)^k,$$
where $k$ is called the length of the segment, denoted by $|s|$.
If ${\mathbf{s}}=\{s_1, \ldots, s_p\}$ is an unordered collection of
segments, define $\varrho(\mathbf{s})$ be the partition associated
with the sequence $(|s_1|,\ldots,|s_p|)$. That is,
$\varrho({\mathbf{s}})=(|s_{i_1}|,\ldots,|s_{i_p}|)$ with
$|s_{i_1}|\geq\ldots\geq|s_{i_p}|$, where $|s_{i_1}|,\ldots,
|s_{i_p}|$ is a permutation of $|s_{1}|,\ldots, |s_{p}|$. We also
call $|{\mathbf{s}}|=\sum_{1\leq i\leq p}|s_i|$ the length of
$\mathbf{s}$.

Let $$C_{r,n}=\{{\mathbf{s}}=\{s_1, \ldots, s_p\}\;|\; p\geq 1,\sum_{1\leq i\leq p}|s_i|=r ,
 |s_i|\leq n, \forall i,
\}.$$
Then
$C_{r,n}=$$\cup_{\lb\in\Lambda^+(n,r) }$$C_{\lambda}$,
where $C_{\lb}=\{{\mathbf{s}}\in C_{r,n}\;|\;
\varrho(\mathbf{s})=\lb\}$. Note that if $\lb=(\lb_1,\ldots, \lb_n)\in \Lambda^+(n,r)$
and there is some $1\leq i\leq n$ such that $\lb_i>n$,
then $C_\lb=\emptyset$.

\begin{lem}[\cite{DDF}, \cite{F}]\label{82}
$C_{r,n} $ is a parameter set of finite
dimensional non-isomorphic simple $\whS(n,r)_\bbc$-modules.

\end{lem}

\begin{pf}

In \cite{F}, finite dimensional simple
$\whS(n,r)_\bbc$-modules are parameterized by the set of dominant polynomials. It is proved
in \cite[Thm.~4.4.2]{DDF} and \cite[Thm.~6.6]{DD} that the set of dominant polynomials and
$C_{r,n}$ coincide.
\end{pf}

For $\lb=(\lb_1,\ldots,\lb_n)\in \Lambda^+(n,r)$, we define its
dual partition $\lb'$ as follows
$$\lb'=(\lb_1',\ldots, \lb_n'),$$where
$\lb_i'=|\{j|\lb_j\geq i\}|$, for $1\leq i\leq n$. Note that for $\lb\in \Lambda^+(n,r)$,
its dual partition $\lb'\in \Lambda^+(n,r)$ if and only if $\lb_i\leq n$ for all $1\leq i\leq n$.

\begin{thm}\label{83} For arbitrary $r, n\in\bbn$, there is a bijection between $C_{r,n}$
and $\Omega_{r,n}$.
\end{thm}

\begin{pf}

Take $\lb=(\lb_1, \ldots, \lb_i, 0,\ldots,0)\in\Lambda^+(n,r) $ with $\lb_i\neq 0$ and $\lb_j\leq n$ for $1\leq j\leq i$.
Then
$${\mathbf{s}}=(\underbrace{a_1, \ldots,
a_1}_{\lb_1}, \ldots, \underbrace{a_i, \ldots,
a_i}_{\lb_i})\in C_{\lb},$$ where $a_t\neq 0$
for $1\leq t\leq i$.

Define
$$\aligned
\phi: C_{r,n}&\longrightarrow \Omega_{r,n}, \\
{\mathbf{s}}=(\underbrace{a_1, \ldots, a_1}_{\lb_1}, \ldots, \underbrace{a_i, \ldots,
a_i}_{\lb_i})&\longmapsto (b_1, \ldots, b_i),
\endaligned$$
 where $b_{i_j}=e_{i_j}$ is the $j$-th elementary symmetric function in $a_{i_{1}},
\ldots, a_{i_s}$ if $\lb_{i_1}=\cdots=\lb_{i_s}$ for $1\leq i_1,\ldots, i_s\leq i$.

It is easy to check that $\phi({\mathbf{s}})\in \Omega_{\lb'}$ and $\phi$ is a bijection between $C_{r,n}$ and $\Omega_{r,n}$. The proof is
completed.

\end{pf}

\begin{example} Let us consider some examples. Let $\lb=(3,2,1)$. Then
$$C_\lb=\{\underline{a}=(a_1, a_1, a_1, a_2, a_2, a_3)\;|\; a_i\in\bbc, \;\text{for}\;1\leq i\leq 3, a_1a_2a_3\neq 0\}.$$
Since $3\neq 2\neq 1$, by Theorem \ref{83},
$$\phi(\underline{a})=(a_1, a_2, a_3)\in \Omega_{\lb'},$$
where $\lb'=(3,2,1)$.

Let $\mu=(3,3,1)$. Then
 $$C_\mu=\{\underline{a'}=(a_1, a_1, a_1, a_2, a_2, a_2, a_3)\;|\;a_i\in\bbc, \;\text{for}\;1\leq i\leq 3, a_1a_2a_3\neq 0\}.$$
Since $3=3\neq 1$, by Theorem \ref{83},
$$\phi(\underline{a'})=(a_1+a_2, a_1a_2, a_3)\in\Omega_{\mu'},$$
where $\mu'=(3,2,2)$.
Let $\nu=(4,2,1)$. Since $\nu_1=4>3$, this implies that $C_\nu=\emptyset$. Its dual partition $\nu'=(3,2,1,1)$ does not
belong to $\Lambda^+(3,7)$.

\end{example}

\end{document}